\documentclass[12pt]{article}
\usepackage{amsmath, enumerate, amssymb, bbm, amssymb,amsthm,bm}

\newtheorem{theorem}{Theorem}[section]
\newtheorem{lemma}[theorem]{Lemma}
 
\newtheorem{remark}[theorem]{Remark}
\newtheorem{proposition}[theorem]{Proposition}
\newtheorem{assumption}[theorem]{Assumption}

\numberwithin{equation}{section}

\def\grad{{\nabla}}
\def\proof{{\medskip\noindent {\bf Proof. }}}
\def\qed{{\hfill $\square$ \bigskip}}
\def\longproof#1{{\medskip\noindent {\bf Proof #1.}}}

 \def\sE {{\cal E}} \def\sF {{\cal F}}
  
  \def\sL {{\cal L}}
 \def\sN {{\cal N}}

 \def\bE {{\Bbb E}}

 \def\bN {{\Bbb N}} 
\def\bP {{\Bbb P}}  \def\bR {{\Bbb R}}

\begin{document}

\title{\bf Heat kernel estimates and Harnack inequalities for some Dirichlet forms with non-local~part}

\author{Mohammud Foondun\footnote{Department of Mathematics, University of Utah, Salt Lake City, UT 08112; mohammud@math.utah.edu}}

\date{}

\maketitle
\abstract{We consider the Dirichlet form given by 
\begin{eqnarray*}
\sE(f,f)&=&\frac{1}{2}\int_{\bR^d}\sum_{i,j=1}^d a_{ij}(x)\frac{\partial f(x)}{\partial x_i} \frac{\partial f(x)}{\partial x_j} dx\\
&+&\int_{\bR^d\times \bR^d} (f(y)-f(x))^2J(x,y)dxdy.
\end{eqnarray*}
Under the assumption that the $\{a_{ij}\}$ are symmetric and uniformly elliptic and with suitable conditions on $J$, the nonlocal part, we obtain upper and lower bounds on the heat kernel of the Dirichlet form.  We also prove a Harnack inequality  and  a regularity theorem for functions that are harmonic with respect to $\sE$.}

\vglue1.0truein
{\small

}
\newpage

\section{Introduction}
The main aim of this article is prove a Harnack inequality and a regularity estimate for harmonic functions with respect to some Dirichlet forms with non-local part.  More precisely, we are going to consider the following Dirichlet form
\begin{eqnarray}\label{defdirich}
\sE(f,f)&=&\frac{1}{2}\int_{\bR^d}\sum_{i,j=1}^d a_{ij}(x)\frac{\partial f(x)}{\partial x_i} \frac{\partial f(x)}{\partial x_j} dx \nonumber\\
&+& \int_{\bR^d}\int_{\bR^d} (f(y)-f(x))^2J(x,y)dxdy,
\end{eqnarray}
where $a_{ij}:\bR^d\rightarrow \bR$ and $J:\bR^d\times \bR^d\rightarrow \bR$ satisfy some suitable assumptions; see Assumptions \ref{ass1} and \ref{ass2} below. The domain $\sF$  of the Dirichlet form $\sE$ is defined as the closure with respect to the metric $\sE_1^{1/2}$ of $C^1$-functions on $\bR^d$ with compact support, where $\sE_1$ is given by: $\sE_1(f,f):=\sE(f,f)+\int_{\bR^d}f(x)^2dx$.

The local part of the above form corresponds to the following elliptic operator
\begin{equation}\label{diff}
\sL=\sum_{i,j=1}^d\frac{\partial}{\partial x_i}\left(a_{ij}(x)\frac{\partial}{\partial x_j}\right)
\end{equation}
which was studied in the papers of E.DeGiorgi\cite{De}, J.Nash\cite{N} and J.Moser\cite{M1,M2} as well as in many others.  They showed that under the assumptions that the matrix $a(x)=(a_{ij}(x))$ is symmetric and uniformly elliptic, harmonic functions with respect to $\sL$ behave much like those with respect to the usual Laplacian operator.  This holds true even though the coefficients $a_{ij}$ are assumed to be measurable only.  The above Dirichlet form given by $(\ref{defdirich})$ has a probabilistic interpretation in that it represents a discontinuous process with the local part representing the continuous part of the process while the non-local part represents the jumps of the process.  We call $J(x,y)$ the jump kernel of the Dirichlet form. It represents the intensity of jumps from $x$ to $y$.

In a way, this paper can be considered as the analogue of our earlier paper \cite{Fo1} where the following operator was considered:
\begin{eqnarray}\label{nondiv}
\sL f(x)&=&\frac{1}{2}\sum_{i,j=1}^d a_{ij}(x)\frac{\partial^2f(x)}{\partial x_i\partial x_j}+\sum_{i=1}^d b_i(x)\frac{\partial f(x)}{\partial x_i}\nonumber\\
&+&\int_{\bR^d\backslash\{0\}}[f(x+h)-f(x)-1_{(|h|\leq1)}h\cdot \grad f(x)]n(x,h)dh.
\end{eqnarray}

In that paper, a Harnack inequality as well as a regularity theorem were proved.  The methods employed were probabilistic and there we related the above operator to a process via the martingale problem of Stroock and Varadhan, whereas here the probabilistic interpretation is given via the theory described in \cite{FOT}.

The study of elliptic operators has a long history.  E. DeGiorgi\cite{De}, J. Nash\cite{N} and J. Moser\cite{M1}, among others, made significant contributions to the understanding of elliptic operators in divergence form.  In \cite{KS1} Krylov and Safonov gave a probabilistic proof of the Harnack inequality as well as a regularity estimate for elliptic operators in non-divergence form.  

While there has been a lot of research concerning differential operators, not much has been done for non-local operators.  It is only recently that Bass and Levin \cite{BL1} proved a Harnack inequality and a continuity estimate  for harmonic functions with respect to some non-local operators.  More precisely, they considered the following operator
\begin{equation}\label{discont}
\sL f(x)=\int_{\bR^d\backslash\{0\}}[f(x+h)-f(x)]\frac{n(x,h)}{|h|^{d+\alpha}}dh,
\end{equation}
where $n(x,h)$ is a strictly positive bounded function satisfying $n(x,h)=n(x,-h)$.  Since then, non-local operators have received considerable attention. For instance in \cite{BK1}, Harnack inequalities were established for variants of the above operator.  Also, Chen and Kumagai \cite{CK} established some heat kernel estimates for stable-like processes in $d$-sets as well as a parabolic Harnack inequality for these processes and in \cite{CK1}, the same authors established heat kernel estimates for jump processes of mixed type in metric spaces. Non-local Dirichlet forms representing pure jump processes have also been recently studied in \cite{BBCK} where bounds for the heat kernel and Harnack inequalities were established.  A special case of the Dirichlet form given by (\ref{defdirich}) was studied by Kassmann in \cite{kas} where a weak Harnack inequality was established. Related work on discontinuous processes include \cite{CS}, \cite{CKS1}, \cite{CKS2}, \cite{SV05}  and \cite{RRV}.

At this point of the introduction it is pertinent to give some more details about the differences between this paper and the results in some related papers.
\begin{itemize}
\item In \cite{kas} a weak Harnack inequality was established and the jump kernel was similar to the one defined in (\ref{discont}) but with index $\alpha \in [1,2)$.  There, the techniques used were purely analytic while here the method used is more probabilistic. This allows us to prove the Harnack inequality and continuity estimate for a much wider class of jump kernels.
\item In \cite{BBCK}, a purely non-local Dirichlet form was considered. The jump kernel considered there satisfies a lower and an upper bound.  Here because of the presence of the local part, no lower bound is required.  The intuitive reason behind this is that since we have a uniformly elliptic local part, the process can move even if there is no jump.  This also agrees with the fact that our results should hold when the jump kernel is identically zero.   
\item A parabolic Harnack inequality was also proved in \cite{BBCK}.  Their result holds on balls with large radius $R$, while here we prove  the Harnack inequality for small $R$ only.  Moreover, in \cite{BBCK} the authors considered processes with small jumps only. Here, our processes are allowed to have big jumps.
\item For our Harnack inequality to hold, we need assumption \ref{ass2}(c) below. This assumption is modeled after the one introduced in \cite{BK1}. Thus with this assumption, our result covers the case when the jump kernel $J(x,y)$ satisfies 
\begin{equation*}
\frac{k_1}{|x-y|^{d+\alpha}}\leq J(x,y)\leq \frac{k_2}{|x-y|^{d+\beta}},\quad{\rm where}\quad0<\alpha<\beta<2,
\end{equation*}
and the $k_i$s are positive constants.  Here, unlike in \cite{BK1}, there is no restriction on $\beta-\alpha$.
\item In a recent preprint \cite{CKK}, Chen, Kim and Kumagai looked at truncated jump processes whose kernel is given by the following
\begin{equation*}
J(x,y)=\frac{c(x,y)}{|x-y|^{d+\alpha}}1_{(|x-y|\leq \kappa)},
\end{equation*}
where $\alpha \in(0,2)$, $\kappa$ is a positive constant and $c(x,y)$ is bounded below and above by positive constants.  The results proved in that paper include sharp heat kernel estimates as well as a parabolic Harnack inequality. The jump kernel studied here includes the ones they study, but since the processes considered here include a continuous part, the results are different.
\end{itemize}
We now give a plan of our article.  In Section 2, we give some preliminaries and state the main results. We give upper and lower bounds for the heat kernel associated to the Dirichlet form in Section~3.  In Section 4, we prove some estimates which will be used in the proof of the regularity theorem and  the Harnack inequality.  In Section 5, a proof of the regularity theorem is given.  A proof of the Harnack inequality is given in Section 6.

\section{Statement of results}
We begin this section with some notations and preliminaries. $B(x,r)$ and $B_r(x)$ will both denote the ball of radius $r$ and center $x$.  The letter $c$ with subscripts will denote positive finite constants whose exact values are unimportant.  The Lebesgue measure of a Borel set A will be denoted by $|A|$.  We consider the Dirichlet form defined by (\ref{defdirich}) and make the following assumptions:
\begin{assumption} \label{ass1}
We assume that the matrix $a(x)=(a_{ij}(x))$ is symmetric and uniformly elliptic.  In other words, there exists a positive constant $\Lambda$ such that the following holds:
\[ \Lambda^{-1}|y|^2\leq \sum_{i,j=1}^dy_ia_{ij}(x)y_j\leq \Lambda|y|^2,\hskip10mm \forall\,\, x,y \in \bR^d.\]
\end{assumption}
We also need the following assumption on the nonlocal part of the Dirichlet form.
\begin{assumption}\label{ass2}
\
\begin{enumerate}[(a)]
\item  There exists a positive function $\tilde{J}$ such that $J(x,y)1_{(|x-y|\leq1)}\leq \tilde{J}(|x-y|)1_{(|x-y|\leq1)}$ for $x,y \in \bR^d$.  Moreover, 
\[\int_{|x-y|\leq 1}|x-y|^2 \tilde{J}(|x-y|)dy\leq K_1 \hskip6mm{\rm and}\hskip6mm \int_{|x-y|>1}J(x,y)dy \leq K_2, \hskip6mm\forall \,\, x\in \bR^d,\]
where  $K_1$ and $K_2$ are positive constants.
\item The function $J(x,y)$ is symmetric, that is, 
\[J(x,y)=J(y,x) \hskip10mm \forall \,\,x,\,y\in \bR^d,\]
\item Let $x_0\in \bR^d$ be arbitrary  and $r\in (0,1]$, then whenever $x,\,y\in B(x_0,r/2)$ and $z\in B(x_0,r)^c$, we have 
\begin{equation*}
J(x,z)\leq k_rJ(y,z),
\end{equation*}
with $k_r$ satisfying $1<k_r\leq \kappa r^{-\beta}$, where $\kappa$ and $\beta$ are constants.
\end{enumerate}
\end{assumption}
In probabilistic terms, $J(x,y)$ can be thought as the intensity of jumps from $x$ to $y$. 
Our method is probabilistic, so we need to work with a process associated with our Dirichlet form.  
The following lemma gives conditions for the existence of a process and its density function.
We say that a Dirichlet form $\sE$ satisfies a \emph{Nash inequality} if 
\begin{equation*}
\|f\|_2^{2(1+\frac{2}{d})} \leq c\sE^{Y^\lambda}(f,f)\|f\|_1^{4/d},
\end{equation*}
where $f\in \sF$ and $c$ is a positive constant.  For an account of various forms of Nash inequalites, see \cite{CKS}.  For a definition of \emph{regular} Dirichlet form, the reader is referred to page 6 of \cite{FOT}. 

\begin{lemma}\label{existence:1}
Suppose that the Dirichlet form is regular and satisfies a Nash inequality, then there exists a process $X$ with a transition density function $p(t,x,y)$ defined on $(0,\infty)\times \bR^d\backslash \sN \times \bR^d\backslash\sN$ satisfying $P(t,x,dy)=p(t,x,y)dy$, where $P(t,x,dy)$ denotes the transition probability of the process $X$ and $\sN$ is a set of capacity zero.
\end{lemma}
\proof{The existence of such a process follows from Theorem 7.2.1 of \cite{FOT} while the existence of the probability density is a consequence of Theorem 3.25 of \cite{CKS}.\qed}

For the rest of the paper, $\sN$ will denote the set of capacity zero, as defined in the above Lemma. 
For any Borel set $A$, let
\[T_A=\inf\{t: X_t\in A\}, \hskip15mm \tau_A=\inf\{t:X_t\notin A \}  \]
be the first hitting time and first exit time, respectively, of $A$.  We say that the function $u$  is harmonic in a domain $D$ if $u(X_{t\wedge \tau_{D}})$ is a $\bP^{x}$-martingale for each $x\in D$.  Since our process is a discontinuous process, we define
\[X_{t-}=\lim_{s\uparrow t}X_s, \hskip10mm{\rm and}\hskip10mm \Delta X_t=X_t-X_{t-}.\]
Here are the main results:
\begin{theorem}\label{theo1}
Suppose Assumptions ~\ref{ass1}, 2.2(a) and 2.2(b) hold. Let $Y$ denote the process associated with the Dirichlet form defined by (\ref{defdirich}) but with jump kernel given by $J(x,y)1_{(|x-y|\leq 1)}$ and null set, $\sN'$. Then there exists a constant $c_1>0$ depending only on $\Lambda$ and the $K_is$ such that for all $x,\,y\in \bR^d\backslash \sN'$ and for all $t\in(0,1]$, the transition density function $p^Y(t,x,y)$ satisfies 
\[p^Y(t,x,y)\leq c_1t^{-\frac{d}{2}}e^{-|x-y|}.\]
\end{theorem}

\begin{theorem}\label{theo2}
Suppose Assumptions ~\ref{ass1}, 2.2(a) and 2.2(b) hold. Let $p(t,x,y)$ denote the transition density function of the process $X$.  Then there exist positive constants $c_1$ and $\theta$ such that
 \[ p(t,x,y) \geq c_1t^{-\frac{d}{2}} \hskip8mm{\rm if} \hskip8mm |x-y|^2\leq \theta t,\hskip8mm {\rm where}\,\,x\,,y\in \bR^d\backslash \sN. \]
\end{theorem}

\begin{theorem}\label{theo3}
Suppose Assumptions ~\ref{ass1}, 2.2(a) and 2.2(b) hold. Let $z_0\in \bR^d$ and $R\in (0,1]$. Suppose u is a function which is bounded in $\bR^d$ and harmonic in $B(z_0,R)$ with respect to the Dirichlet form $(\sE,\sF)$.  Then there exists $\alpha \in (0,1)$ and $C>0$ depending only on $\Lambda$ and the $K_is$ such that 
\[|u(x)-u(y)|\leq C\|u\|_\infty \left(\frac{|x-y|}{R}\right)^\alpha, \hskip10mm x,\,y\in B(z_0,R/2)\backslash \sN. \]
\end{theorem}

\begin{theorem}\label{theo4}
Suppose Assumptions 2.1 and 2.2 hold.  Let $z_0\in \bR^d$ and $R\in (0,1]$.  Suppose $u$ is nonnegative and bounded on $\bR^d$ and harmonic with respect to the Dirichlet form $(\sE,\sF)$ on $B(z_0,R)$.  Then there exists $C>0$ depending only on $\Lambda$, $\kappa$, $\beta$,  $R$ and the $K_is$ but not on $z_0$, $u$ or $\|u\|_\infty$ such that 
\[ u(x)\leq Cu(y), \hskip 10mm x,\,y\in B(z_0,R/2)\backslash \sN. \]
\end{theorem}

We mention that the main ideas used for the proof of the above theorem appear in \cite{BL1}. Note that Assumption \ref{ass2}(c) is crucial for the Harnack inequality to hold.  In fact, an example in the same spirit as that in \cite{Fo1} can be constructed so that the Harnack inequality fails for a Dirichlet form with a jump kernel not satisfying Assumption \ref{ass2}(c).  We do not reproduce this example here because the only difference is that here, we require the process to be symmetric while in \cite{Fo1}, the process is not assumed to be symmetric.

We make a few more comments about some of the assumptions in the above theorem.  We require that the local part is uniformly elliptic and as far as we know, our method does not allow us to relax this condition.  Moreover, as shown in \cite{kas1}, the nonnegativity assumption cannot be dropped. In that paper, the author constructs an example (violating the nonnegativity assumption) which shows that the Harnack inequality can fail for non-local operators.


\section{Upper and lower bounds for the heat kernel}

The main goal of this section is to prove some upper and lower bounds on the heat kernel.  The upper bound on the heat kernel estimate follows from a Nash inequality which is proved in Proposition \ref{generalupper}. For more information about the relation between Nash inequalities and heat kernel estimates, see \cite{CKS}. As for the lower bound, we use Nash's original ideas, see \cite{N}.  Since we are dealing with operators which are not local, we also need some ideas which first appeared in \cite{BBCK}. The paper \cite{SS} also contain some useful information on how to deal with local operators.

We start off this section by proving the regularity of the Dirichlet form $(\sE,\sF)$.  Let $H^1(\bR^d)$ denote the Sobolev space of order $(1,2)$ on $\bR^d$.  In other words, $H^1(\bR^d):=\{f\in L^2(\bR^d): \grad f\in L^2(\bR^d) \}.$

\begin{proposition}\label{regular}
Let $(\sE,\sF)$ be defined by (\ref{defdirich}) . Then,
\[\sF=H^1(\bR^d)=\{f\in L^2(\bR^d): \grad f\in L^2(\bR^d) \}.\]
\end{proposition}\label{regularity}
\proof{We assume that $f$ is continuous with compact support, $K\subset \bR^d$.  Let us write 
\[ \sE(f,f)=\sE_c(f,f)+\sE_d(f,f),\]
where
\begin{eqnarray*}
\sE_c(f,f)&=&\frac{1}{2}\int_{\bR^d}\sum_{i,j=1}^d a_{ij}(x)\frac{\partial f(x)}{\partial x_i} \frac{\partial f(x)}{\partial x_j} dx, \nonumber\\
\sE_d(f,f)&=&\int_{\bR^d\times \bR^d} (f(y)-f(x))^2J(x,y)dxdy.
\end{eqnarray*}
From Assumption \ref{ass1}, we see that if $\grad f\in L^2(\bR^d)$ then $\sE_c(f,f)\leq c_1\|\grad f\|_2^2$.  As for the discontinuous part, we have
\begin{eqnarray*}
\sE_d(f,f)&\leq&\iint_{(B(R)\times B(r))\cap(|x-y|\leq 1)} (f(y)-f(x))^2J(x,y)dxdy\\
&+&\iint_{(B(R)\times B(r))\cap(|x-y|>1)} (f(y)-f(x))^2J(x,y)dxdy\\
&+&2\iint_{B(R)^c\times B(r)} (f(y))^2J(x,y)dxdy\\
&=& I_1+I_2+I_3,
\end{eqnarray*}
where $ B(r)$ and $ B(R)$ are balls with a common center but with radius $r$ and $R$ respectively, satisfying $K\subset B(r) \subset B(R)$ and $R-r>1$.  We consider the term $I_1$ first. Recall that from Assumption \ref{ass2}(a), we have
\begin{equation}\label{eq:levy}
I_1\leq \iint_{(B(R)\times B(r))\cap(|x-y|\leq1)} (f(y)-f(x))^2 \tilde{J}(|x-y|)dxdy.
\end{equation}
Since the measure $\tilde{J}(|h|)1_{(|h|\leq 1)}dh$ is a L\'evy measure, we can use the L\'evy Khintchine formula(see (1.4.21) of \cite{FOT}) to estimate the characteristic function $\psi$ of the corresponding process as follows
\begin{eqnarray*}
\psi(u)&=&\int_{(|h|\leq 1)}(1-\cos(u\cdot h)) \tilde{J}(|h|)\,dh\\
&\leq&c_1\int_{(|h|\leq 1)}|u|^2|h|^2 \tilde{J}(|h|)\,dh\\
&\leq&c_2|u|^2.
\end{eqnarray*}
We now use a simple substitution, Plancherel's theorem as well as the above inequality to obtain 
\begin{eqnarray*}
I_1&\leq&\iint(f(x+h)-f(x))^2 \tilde{J}(|h|)1_{(|h|<1)}\,dh\,dx\\
&\leq&c_3\int|\hat{f}(u)|^2\psi(u)du\\
&\leq&c_4\int|\hat{f}(u)|^2|u|^2du\\
&\leq&c_4\int|\hat{f}(u)|^2(1+|u|^2)du=c_5(\|f\|^2_2+\|\grad f\|_2^2).
\end{eqnarray*}
In the above $\hat{f}$ denotes the Fourier transform of $f$.  A similar argument is used in the proof of (1.4.24) in \cite{FOT}. 
As for the second term $I_2$, we have
\begin{eqnarray*}
I_2&\leq&4\iint_{(|x-y|>1)} |f(y)|^2J(x,y)dxdy\\
&\leq&c_6\|f\|_2^2.
\end{eqnarray*}
The third term $I_3$ is bounded similarly, that is we have $I_3\leq c_7\|f\|_2^2$.
From the above, we see that if $f\in \{f\in L^2(\bR^d): \grad f\in L^2(\bR^d) \}$, then $\sE_1(f,f)<\infty$. Using uniform ellipticity, we can also conclude that if $\sE_1(f,f)<\infty$, then $f\in \{f\in L^2(\bR^d): \grad f\in L^2(\bR^d) \}$.   We now show that for any $u\in H^1(\bR^d)$, there is a sequence $\{u_n\}\subset C^1(\bR^d)$ such that $u_n\rightarrow u$ in the metric $\sE_1^{1/2}$. Denote by $H_0^1(\bR^d)$, the closure of $C_0^\infty(\bR^d)$ in $H^1(\bR^d)$. Then from Proposition 1.1 on page 210 of \cite{CW}, we have $H^1(\bR^d)=H_0^1(\bR^d)$.  Therefore there exists a sequence of $u_n\in C_0^\infty(\bR^d)\subset C_0^1(\bR^d)$ such that $u_n\rightarrow u$ in $H^1(\bR^d)$.  From the calculations above, we have 
\begin{equation}\label{R1}
\sE(u-u_n,u-u_n)\leq c_8\|\grad (u-u_n)\|^2_2+c_9\|u-u_n\|_2^2.
\end{equation}
Letting $n\rightarrow \infty$, we thus have $\sE(u-u_n,u-u_n)\rightarrow 0$. Thus $u_n$ is $\sE_1$-convergent to $u\in \sF$. This shows that $C^1(\bR^d)$ is dense in $(\sE_1,H^1(\bR^d))$, hence concluding the proof.
}
\qed
\begin{remark} \label{remark1}{\rm
In Chapter 7 of \cite{FOT}, it is shown that for any regular Dirichlet form, there exists a Hunt process whose Dirichlet form is the given regular one.  More precisely, there exists $\sN \subset \bR^d$ having zero capacity with respect to the Dirichlet form $(\sE,\sF)$ and there exists a Hunt process $(\bP^x, X)$ with state space $\bR^d\backslash \sN$.  Moreover, the process is uniquely determined on $\sN^c$.  In other words, if there exist two Hunt processes for which the corresponding Dirichlet forms coincide, then there exist a common proper exceptional set $\sN$ so that the transition functions coincide on $\sN^c$.}
\end{remark}

\begin{remark}\label{remark2}{\rm
We will repeatedly use the following construction due to Meyer(\cite{Me}); see also \cite{BBCK} and \cite{BGK}.  This will enable us to restrict our attention to the process with small jumps only and then incorporate the big jumps later. Suppose that we have two jump kernels $J_0(x,z)$ and $J(x,z)$ with $J_0(x,z)\leq J(x,z)$ and such that for all $x\in \bR^d$,
\[N(x)=\int_{\bR^d}(J(x,z)-J_0(x,z))dz \leq c,\]
where $c$ is a constant.

Let $\sE$ and $\sE_0$ be the Dirichlet forms corresponding to the kernels $J(x,z)$ and $J_0(x,z)$ respectively. If $\overline X_t$ is the process corresponding to the Dirichlet form $\sE_0$, then we can construct a process $X_t$ corresponding to the Dirichlet form $\sE$ as follows.  Let $S_1$ be an exponential random variable of parameter 1 independent of $\overline X_t$, let $C_t=\int_0^tN(\overline X_s)ds$, and let $U_1$ be the first time that $C_t$ exceeds~$S_1$.

At the time $U_1$, we introduce a jump from $X_{U_1-}$ to $y$, where $y$ is chosen at random according to the following distribution:
\[ \frac{J(\overline X_{U_1-},z)-J_0(\overline X_{U_1-},z)}{N(\overline X_{U_1-})}dz.\]

This procedure is repeated using an independent exponential variable $S_2$.  And since $N(x)$ is finite, for any finite time interval we have introduced only a finite number of jumps.  Using \cite{Me}, it can be seen that the new process corresponds to the Dirchlet form $\sE$.  And if $\sN_0$ is the set of zero capacity corresponding to the Dirichlet form $\sE_0$, then $\sN\subset \sN_0$.}

\end{remark}


\subsection{Upper bounds}
Let $Y^\lambda$ be the process associated with the following Dirichlet form:
\begin{eqnarray}\label{D1}
\sE^{Y^\lambda}(f,f)&=&\frac{1}{2}\int_{\bR^d}\sum_{i,j=1}^d a_{ij}(x)\frac{\partial f(x)}{\partial x_i} \frac{\partial f(x)}{\partial x_j} dx \nonumber\\
&+& \iint_{|x-y|\leq \lambda} (f(y)-f(x))^2J(x,y)dxdy,
\end{eqnarray}
so that $Y^{\lambda}$ has jumps of size less than $\lambda$ only.  Let $\sN(\lambda)$ be the exceptional set corresponding to the Dirichlet form defined by (\ref{D1}).  Let $P_t^{Y^\lambda}$ be the semigroup associated with $\sE^{Y^\lambda}$.  We will use the arguments in \cite{FOT} and \cite{CKS}  as indicated in the proof of Lemma \ref{existence:1} to obtain the existence of the heat kernel $p^{Y^\lambda}(t,x,y)$ as well as some upper bounds.
For any $v, \psi\in \sF$, we can define
\begin{eqnarray*}
\Gamma_\lambda[v](x)&=& \frac{1}{2}\grad v\cdot a\grad v+\int_{|x-y|\leq \lambda}(v(x)-v(y))^2J(x,y)dy,\\
D_\lambda(\psi)^2&=&\|e^{-2\psi}\Gamma_\lambda[e^\psi]\|_\infty \vee\|e^{2\psi}\Gamma_\lambda[e^{-\psi}]\|_\infty,
\end{eqnarray*}
and provided that $D_\lambda(\psi) <\infty$, we set 
\begin{eqnarray*}
E_\lambda(t,x,y)&=&\sup\{|\psi(y)-\psi(x)|-tD_\lambda(\psi)^2; D_\lambda(\psi) <\infty\}.
\end{eqnarray*}

\begin{proposition}\label{generalupper}
There exists a constant $c_1$ such that the following holds.
\[p^{Y^\lambda}(t,x,y) \leq c_1t^{-d/2}\exp[-E_{\lambda}(2t,x,y)], \hskip5mm\forall \,\, x,\,y\in \bR^d\backslash \sN(\lambda), \hskip 5mm {\rm and}\,\,\,t\in (0,\infty),\]
where $p^{Y^\lambda}(t,x,y)$ is the transition density function for the process $Y^\lambda$ associated with the Dirichlet form $\sE^{Y^\lambda}$.
\end{proposition}
\proof{
Similarly to Proposition \ref{regular}, we write
\[\sE^{Y^\lambda}(f,f)=\sE^{Y^\lambda}_c(f,f)+\sE^{Y^\lambda}_d(f,f),\]
Since $J(x,y) \geq 0$ for all $x,\,y\in \bR^d$, we have 
\begin{equation}\label{mono}
\sE^{Y^\lambda}(f,f)\geq \sE^{Y^\lambda}_c(f,f).
\end{equation}
We have the following Nash inequality; see Section VII.2 of \cite{Ba3}:
\[ \|f\|_2^{2(1+\frac{2}{d})} \leq c_2\sE^{Y^\lambda}_c(f,f)\|f\|_1^{4/d}.\]
This, together with (\ref{mono}) yields
\[ \|f\|_2^{2(1+\frac{2}{d})} \leq c_2\sE^{Y^\lambda}(f,f)\|f\|_1^{4/d}.\]
Now applying Theorem 3.25 from \cite{CKS}, we get the required result.  
\qed
}

We now estimate $E_\lambda(t,x,y)$ to obtain our first main result.

 \longproof{of Theorem \ref{theo1}}
 Let us write $\Gamma_\lambda$ as 
\[\Gamma_\lambda[v]=\Gamma_\lambda^c[v]+\Gamma_\lambda^d[v],\]
where 
\[\Gamma_\lambda^c[v]= \frac{1}{2}\grad v\cdot a\grad v,\]
and
\[\begin{displaystyle} \Gamma _\lambda^d[v]=\int_{|x-y|\leq \lambda}(v(x)-v(y))^2J(x,y)dy\end{displaystyle}.\]
\\
Fix  $(x_0,y_0)\in \left(\bR^d\backslash \sN(\lambda)\right)\times\left(\bR^d\backslash \sN(\lambda)\right)$. Let $\mu>0$ be constant to be chosen later.  Choose $\psi(x)\in \sF$ such that $|\psi(x)-\psi(y)|\leq \mu|x-y|$ for all $x,y \in \bR^d$. We therefore have the following:
\begin{eqnarray*}
\left|e^{-2\psi(x)}\Gamma_\lambda^d[e^{\psi}](x)\right|&=&e^{-2\psi(x)}\int_{|x-y|\leq \lambda}(e^{\psi(x)}-e^{\psi(y)})^2J(x,y)dy \\
&=&\int_{|x-y|\leq \lambda}( e^{\psi(y)-\psi(x)}-1)^2J(x,y)dy\\
&\leq&c_1\int_{|x-y|\leq \lambda}|\psi(x)-\psi(y)|^2e^{2|\psi(x)-\psi(y)|}J(x,y)dy\\
&\leq&c_1\mu^2e^{2\mu\lambda}\int_{|x-y|\leq \lambda}|x-y|^2J(x,y)dy\\
&=&c_1\mu^2K(\lambda)e^{2\mu\lambda},
\end{eqnarray*}
where $\begin{displaystyle}K(\lambda)=\sup_{x\in \bR^d}\int_{|x-y|\leq \lambda}|x-y|^2J(x,y)dy\end{displaystyle}$.  Some calculus together with the ellipticity condition yields:
\begin{eqnarray*}
\left|e^{-2\psi(x)}\Gamma_\lambda^c[e^{\psi}](x)\right|&=&\frac{1}{2}\left| e^{-2\psi(x)}\grad(e^{\psi(x)})\cdot a\grad(e^{\psi(x)})\right|\\
&=&\frac{1}{2}\left|\grad \psi(x) \cdot a\grad \psi(x)\right| \\
&\leq&\frac{1}{2}\Lambda \|\grad \psi\|_\infty^2\\
&\leq&2\mu^2 \Lambda.
\end{eqnarray*}
Combining the above we obtain
\[\left|e^{-2\psi(x)}\Gamma_\lambda[e^{\psi}](x)\right|\leq c_1\mu^2K(\lambda)e^{2\mu\lambda}+2\mu^2 \Lambda.\]
Since we have similar bounds for $\left|e^{2\psi(x)}\Gamma_\lambda[e^{-\psi}](x)\right|$, we have 
\begin{eqnarray}\label{Im}
-E_\lambda(2t;x,y)&\leq&2tD_\lambda(\psi)^2-|\psi(y)-\psi(x)|\nonumber\\
&\leq&2t\mu^2\left(c_1K(\lambda)e^{2\mu\lambda}+2 \Lambda\right)-\left| \frac{\mu(x-y)\cdot(x_0-y_0)}{|x_0-y_0|}\right|.
\end{eqnarray}
 Taking $x=x_0$, $y=y_0$  and $\mu=\lambda=1$ in the above and using Proposition \ref{generalupper} together with the fact that $t\leq 1$, we obtain
 \[p^Y(t,x_0,y_0)\leq c_2t^{-\frac{d}{2}}e^{-|x_0-y_0|},\]
Since $x_0$ and $y_0$ were taken arbitrarily, we obtain the required result.
\qed

The following is a consequence of Proposition \ref{generalupper} and an application of Meyer's construction.  

\begin{proposition}\label{generaltightness}
Let $r\in(0,1]$.  Then for $x\in \bR^d\backslash \sN$,
\[\bP^x(\sup_{s\leq {t_0r^2}}|X_s-x|>r)\leq \frac{1}{2},\]
where $t_0$ is a small constant.
\end{proposition}

\proof{ The proof is a follow up of that of the Theorem \ref{theo1}, so we refer the reader to some of the notations there. Let $\lambda$ be a small positive constant to be chosen later.  Let $Y^\lambda$ be the subprocess of X having jumps of size less or equal to $\lambda$.   Let $\sE^{Y^\lambda}$ and $p^{Y^\lambda}(t,x,y)$ be the corresponding Dirichlet form and probability density function respectively.  According to Proposition \ref{generalupper}, we have 
 \begin{equation}\label{GT0}
 p^{Y^\lambda}(t,x,y) \leq c_1t^{-d/2}\exp[-E_{\lambda}(2t,x,y)]. 
 \end{equation}
 Taking $x=x_0$ and $y=y_0$ in ($\ref{Im}$) yields
 \begin{equation}
 -E_\lambda(2t;x_0,y_0)\leq 2t\mu\Lambda+2c_2t\mu^2K(\lambda)e^{2\mu\lambda}-\mu|x_0-y_0|
 \end{equation}
 Taking $\lambda$ small enough so that $K(\lambda)\leq \frac{1}{2c_2}$, the above reduces to 
 \begin{eqnarray*}
 -E_\lambda(2t;x_0,y_0)&\leq&2t\mu^2\Lambda+(t/\lambda^2)(\mu\lambda)^2e^{2\mu\lambda}-\mu|x_0-y_0|\\
 &\leq&2t\mu^2\Lambda+(t/\lambda^2)e^{3\mu\lambda}-\mu|x_0-y_0|.
 \end{eqnarray*}
 Upon setting $\mu=\frac{1}{3\lambda}\log\left(\frac{1}{t^{1/2}}\right)$ and choosing $t$ such that $t^{1/2}\leq \lambda^2$, we obtain
 \begin{eqnarray*}
 -E_\lambda(2t;x_0,y_0)&\leq&c_3t^{1/2}(\log t)^2+\frac{t}{\lambda^2}\frac{1}{t^{1/2}}-\frac{|x_0-y_0|}{3\lambda}\log\left(\frac{1}{t^{1/2}} \right)\\
 &\leq&c_3t^{1/2}(\log t)^2+1+\log[t^{|x_0-y_0|/6\lambda}].
 \end{eqnarray*}
 Applying the above to $(\ref{GT0})$ and simplifying
\begin{eqnarray*}
 p^{Y^\lambda}(t,x_0,y_0)&\leq& c_4e^{c_3t^{1/2}(\log t)^2}t^{|x_0-y_0|/6\lambda}t^{-d/2}\\
 &=&c_4e^{c_3t^{1/2}(\log t)^2}t^{|x_0-y_0|/12\lambda-d/2}t^{|x_0-y_0|/12\lambda}\\
 &=&c_4e^{c_3t^{1/2}(\log t)^2}t^{|x_0-y_0|/12\lambda-d/2}e^{\frac{|x_0-y_0|}{12\lambda}\log t}.
  \end{eqnarray*}
For small $t$, the above reduces to 
\begin{equation}\label{GT1}
p^{Y^\lambda}(t,x_0,y_0)\leq c_5t^{|x_0-y_0|/12\lambda-d/2}e^{-c_6|x_0-y_0|/12\lambda}
\end{equation}

Let us choose $\lambda=c_7r/d$ with $c_7<1/24$ so that for $|x_0-y_0|>r/2$, we have  $|x_0-y_0|/12\lambda-d/2>d/2$. Since $t$ is small(less than one), we obtain

\begin{eqnarray*}
\bP^{x_0}(|Y^\lambda_t-x_0|>r/2)&\leq&\int_{|x_0-y|>r/2}c_5t^{|x_0-y|/12\lambda-d/2}e^{-c_3|x_0-y|/12\lambda}dy\\
&\leq&c_5t^{d/2}\int_{|x_0-y|>r/2}e^{-c_3|x_0-y|/12\lambda}dy.
\end{eqnarray*}
We bound the integral on the right hand side to obtain 
\begin{equation*}
\bP^{x_0}(|Y^\lambda_t-x_0|>r/2)\leq c_7t^{d/2}e^{-c_8r}.
\end{equation*}
Therefore there exists $t_1>0$ small enough such that for $0\leq t \leq t_1$, we have 
\begin{equation*}
\bP^{x_0}(|Y^\lambda_t-x_0|>r/2)\leq \frac{1}{8}.
\end{equation*}
We now apply Lemma 3.8 of $\cite{BBCK}$ to obtain
\begin{equation}\label{GT2}
\bP^x(\sup_{s\leq t_1}|Y_s^\lambda-Y_0^\lambda|\geq r)\leq \frac{1}{4} \hskip10mm \forall \,s\in(0,t_1].
\end{equation} 
We can now use Meyer's argument(Remark \ref{remark2}) to recover the process $X$ from $Y^\lambda$. Recall that in our case $J_0(x,y)=J(x,y)1_{(|x-y|\leq \lambda)}$ so that after using  Assumptions \ref{ass2}(a) and choosing $c_7$ smaller if necessary, we obtain
\begin{eqnarray*}
\sup_x N(x)\leq c_9r^{-2},
\end{eqnarray*}
where $c_9$ depends on the $K_is$ and 
\[N(x)=\int_{\bR^d}(J(x,z)-J_0(x,z))dz.\]
Set $t_2=t_0r^2$ with $t_0$ small enough so that $t_2\leq t_1$.  Recall that $U_1$ is the first time at which we introduce the big jump.  We thus have  
\begin{eqnarray*}
\bP^{x_0}(\sup_{s\leq t_2}|X_s-x_0|\geq r)&\leq&\bP^{x_0}(\sup_{s\leq t_2}|X_s-x_0|\geq r, U_1>t_2)+\bP^{x_0}(\sup_{s\leq t_2}|X_s-x_0|\geq r, U_1\leq t_2)\\
&\leq&\bP^{x_0}(\sup_{s\leq t_2}|Y_s^\lambda-x|\geq r)+\bP^{x_0}(U_1\leq t_2)\\
&=& \frac{1}{4}+1-e^{-(\sup N)t_2}\\
&=& \frac{1}{4}+1-e^{-c_9t_0}.
\end{eqnarray*}
By choosing $t_0$ smaller if necessary, we get the desired result.
\qed
}

\begin{remark}
It can be shown that the process $Y^\lambda$ is conservative.  This fact has been used above through Lemma 3.8 of \cite{BBCK}. 
\end{remark}

\subsection{Lower bounds}
The main aim of this subsection is to prove Theorem $\ref{theo2}$.  We are going to use Nash's original ideas as used in \cite{BBCK}, \cite{CKK} and \cite{SS}.  
Let $x_0\in \bR^d$ and $R>0$.  Set
\begin{equation}\label{cutoff}
\phi_R(x)=((1-\frac{|x-x_0|}{R})^+)^2\quad{\rm for\,\,\,all}\quad{x\in \bR^d},
\end{equation}
and recall that 
\begin{eqnarray}\label{eq:bi}
\sE(f,g)&=&\frac{1}{2}\int_{\bR^d}\sum_{i,j=1}^d a_{ij}(x)\frac{\partial f(x)}{\partial x_i} \frac{\partial g(x)}{\partial x_j} dx\nonumber\\
&+&\int_{\bR^d\times \bR^d} (f(y)-f(x))(g(y)-g(x))J(x,y)dxdy,
\end{eqnarray}
for $f,g\in \sF$.
We begin with the following technical result.
\begin{proposition}\label{prop:tech}
\begin{enumerate}[(a)]
\item There exists a positive constant $c_1$ such that 
$\begin{displaystyle}\left|\frac{\partial p(t,x,y)}{\partial t}\right|\leq c_1t^{-1-d/2}\end{displaystyle}$ for all $t>0$,
\item Fix  $y_0 \in \bR^d\backslash \sN$ and $\epsilon>0$. If $F(t)=\int \phi_R(x)\log p_{\epsilon}(t,x,y_0)dx$, then 
\begin{equation}\label{t1}
F'(t)=-\sE\left(p(t,\cdot, y_0), \frac{\phi_{R}(\cdot)}{p_{\epsilon}(t,\cdot, y_0)}  \right),
\end{equation}
where $p_{\epsilon}(t,x,y):=p(t,x,y)+\epsilon$.
\end{enumerate}
\end{proposition}
\proof{The proof of the first part of the proposition is omitted because it is similar to the proof of Lemma 4.1 of \cite{BBCK}.  We now give a proof of the second part.  We first need to argue that the right hand side of \eqref{t1} makes sense.  The second step is to show the equality \eqref{t1}.\\
{\it Step 1}: By Proposition \ref{regular}, it suffices to show that $\frac{\phi_{R}(\cdot)}{p_{\epsilon}(t,\cdot, y_0)}\in L^2(\bR^d)$ and $\grad \left(\frac{\phi_{R}(\cdot)}{p_{\epsilon}(t,\cdot, y_0)}\right) \in L^2(\bR^d)$. $\frac{\phi_{R}(\cdot)}{p_{\epsilon}(t,\cdot, y_0)}\in L^2(\bR^d)$ follows from the definition of $\phi_R(\cdot)$ and the fact that $p_{\epsilon}(t,x,y)$ is strictly positive.  By Lemma 1.3.3 of \cite{FOT}, we have that $p(t, \cdot, y_0)\in \sF$. Using some calculus, we can write 
\begin{equation*}
\grad \left(\frac{\phi_{R}(\cdot)}{p_{\epsilon}(t,\cdot, y_0)}\right)=\frac{p_{\epsilon}(t,\cdot, y_0)\grad \phi_R(\cdot)-\phi_R(\cdot)\grad p_{\epsilon}(t,\cdot, y_0)}{(p_\epsilon(t,\cdot, y_0))^2}.
\end{equation*}
The above display together with the fact that $p(t, \cdot, y_0)\in \sF$ and the positivity of $p_{\epsilon}(t,x,y)$ show that  $\grad \left(\frac{\phi_{R}(\cdot)}{p_{\epsilon}(t,\cdot, y_0)}\right) \in L^2(\bR^d)$.\\
{\it Step 2}: We write $(f,g)$ for $\int f(x)g(x)dx$.  By Lemma 1.3.4 of \cite{FOT}, we have
\begin{eqnarray*}
-\sE\left(p(t,\cdot, y_0), \frac{\phi_{R}(\cdot)}{p_{\epsilon}(t,\cdot, y_0)}  \right)&=&\lim_{h\rightarrow 0}\frac{1}{h}\left(p(t+h,\cdot,y_0)-p(t,\cdot, y_0),\frac{\phi_{R}(\cdot)}{p_{\epsilon}(t,\cdot, y_0)}\right)\\
&=&\lim_{h\rightarrow 0}\frac{1}{h}\left(p_{\epsilon}(t+h,\cdot,y_0)-p_{\epsilon}(t,\cdot, y_0),\frac{\phi_{R}(\cdot)}{p_{\epsilon}(t,\cdot, y_0)}\right)\\
&=&\lim_{h\rightarrow 0}\frac{1}{h}\int \phi_R(x)\left(\frac{p_{\epsilon}(t+h,x,y_0)}{p_{\epsilon}(t,x, y_0)}-1\right)dx.
\end{eqnarray*}
Taking into consideration the upper and lower bounds on $p_{\epsilon}(t,x,y)$, we see that the right hand side of the above is well defined. 
We have 
\begin{equation*}
F'(t)=\lim_{h\rightarrow 0}\frac{1}{h}\int(\log p_{\epsilon}(t+h,x,y_0)-\log p_{\epsilon}(t,x,y_0))\phi_R(x)dx.
\end{equation*}
Set 
\begin{equation*}
D(h)=[\log p_{\epsilon}(t+h,x,y_0)-\log p_{\epsilon}(t,x,y_0)-(\frac{p_{\epsilon}(t+h,x,y_0)}{p_{\epsilon}(t,x,y_0)}-1)]\phi_R(x).
\end{equation*}
This gives
\begin{equation*}
D'(h)=\frac{\partial}{\partial t}p_{\epsilon}(t+h,x,y_0)(p_{\epsilon}(t,x,y_0)-p_{\epsilon}(t+h,x,y_0))\frac{\phi_R(x)}{p_{\epsilon}(t+h,x,y_0)p_{\epsilon}(t,x,y_0)}.
\end{equation*}
Using the mean value theorem, $D(h)/ h=D'(h^\ast)$ where $h^\ast=h^\ast(x,y_0,h)\in (0,h)$. The bounds on $p_{\epsilon}(t,x,y)$ imply that $D(h)/ h$ tends to $0$ for $x \in B_R(x_0)$ as $h\rightarrow 0$. An application of the  dominated convergence theorem then yields the desired result.
\qed}

We will need the following Poincar\'e inequality.  A proof can be found in \cite{SS}.

\begin{proposition}
Consider the function defined by $(\ref{cutoff})$, there exists a constant $c_1$ not depending on $R$, $f$ and $y_0$, such that 
\begin{equation}\label{poincare}
\int_{\bR^d}|f(x)-\overline{f}|^2\phi_R(x)dx \leq c_1R^2\int_{\bR^d} |\grad f(x)|^2\phi_R(x)dx, \hskip 10mm f\in C_b^\infty(\bR^d).
\end{equation}
where 
\[ \overline{f}=\int_{\bR^d}f(x)\phi_R(x)dx/\int \phi_R(x)dx.\]
\end{proposition}
\longproof{of Theorem \ref{theo2}}: Let $R>0$  and take an arbritary $\epsilon >0$.  Fix $z\in \bR^d$ such that $z\in B_R(0)$ and define $\phi_R(x)=((1-|x|/R)^+)^2$ for $x\in \bR^d$. Set
\begin{eqnarray*}
p_\epsilon(t,x,y)&=&p(t,x,y)+\epsilon,\\
u(t,x)&=&|B_R(0)|p(tR^2,z,x),\\
u_{\epsilon}(t,x)&=&|B_R(0)|p_{\epsilon}(tR^2,z,x),\\
r_{\epsilon}(t,x)&=&\frac{u_{\epsilon}(t,x)}{\phi_{R}(x)^{1/2}},\\
\mu(R)&=&\int\phi_{R}(x)dx,\\
G_\epsilon(t)&=&\mu(R)^{-1}\int\phi_{R}(x)\log u_{\epsilon}(t,x)dx.
\end{eqnarray*}
Using part(b) of Proposition \ref{prop:tech}, we then have
\begin{eqnarray}\label{G1}
\mu(R)G'_\epsilon(t)&=&-R^2\sE\left(u(t,\cdot), \frac{\phi_{R}(\cdot)}{u_{\epsilon}(t,\cdot)}\right)\nonumber\\
 &=&-R^2\left[\sE^{c}\left(u(t,\cdot), \frac{\phi_{R}(\cdot)}{u_{\epsilon}(t,\cdot)}\right)+\sE^{d}\left(u(t,\cdot), \frac{\phi_{R}(\cdot)}{u_{\epsilon}(t,\cdot)}\right)\right]\nonumber\\
 &=&-R^2(I_1+I_2),
  \end{eqnarray}
 where $\sE^{c}$ and $\sE^{d}$ are the local and non-local parts of the Dirichlet form $\sE$ respectively.  Let us look at $I_2$ first. By considering the local part of \eqref{eq:bi} and doing some algebra, we obtain 
\begin{eqnarray*}
I_2&=&\iint\frac{[u(t,y)-u(t,x)]}{u_{\epsilon}(t,x)u_{\epsilon}(t,y)}[u_{\epsilon}(t,x)\phi_R(y)-u_{\epsilon}(t,y)\phi_R(x)]J(x,y)dxdy\\
&=&\iint\frac{[u_{\epsilon}(t,y)-u_{\epsilon}(t,x)]}{u_{\epsilon}(t,x)u_{\epsilon}(t,y)}[u_{\epsilon}(t,x)\phi_R(y)-u_{\epsilon}(t,y)\phi_R(x)]J(x,y)dxdy.\\
\end{eqnarray*}
Note that for $A>0$, the following inequality holds
\begin{equation}\label{ineq:log}
A+\frac{1}{A}-2\geq (\log A)^2.
\end{equation}
We now set $a=u_\epsilon(t,y)/u_\epsilon(t,x)$ and $b=\phi_R(y)/\phi_R(x)$ and observe that
\begin{equation*}
\frac{[u_{\epsilon}(t,y)-u_{\epsilon}(t,x)]}{u_{\epsilon}(t,x)u_{\epsilon}(t,y)}[u_{\epsilon}(t,x)\phi_R(y)-u_{\epsilon}(t,y)\phi_R(x)]
\end{equation*}
\begin{eqnarray*}
\phantom{dddddddddddddddd}&=&\phi_R(x)[b-\frac{b}{a}-a+1]\\
&=&\phi_R(x)[(1-b^{1/2})^2-b^{1/2}(\frac{a}{b^{1/2}}+\frac{b^{1/2}}{a}-2)].
\end{eqnarray*}
Applying inequality \eqref{ineq:log} with $A=a/\sqrt{b}$ to the above equality, we obtain 
\begin{equation*}
I_2\leq \iint[(\phi_{R}(x)^{\frac{1}{2}}-\phi_{R}(y)^{\frac{1}{2}})^2-(\phi_{R}(x)\wedge\phi_{R}(y))\left(\log \frac{r_{\epsilon}(t,y)}{r_{\epsilon}(t,x)}\right)^2]J(x,y)dxdy.
\end{equation*}
See Proposition 4.9 of \cite{BBCK} where a similar argument is used.  We also have 
\[\iint \left[(\phi_{R}(x)\wedge\phi_{R}(y))(\log (r_{\epsilon}(t,y)/ r_{\epsilon}(t,x)))^2J(x,y)dxdy\right]\geq 0.\]
Assumption (\ref{ass2})(a) and the definition of $\phi_{R}(x)$  give the following
\[\iint(\phi_{R}(x)^{1/2}-\phi_{R}(y)^{1/2})^2J(x,y)dxdy \leq c_3|B_R(0)|/{R^2}.\]
Hence we have $I_2\leq c_3|B_R(0)|/{R^2}$. As for the continuous part $I_1$, we use some calculus to obtain 
\begin{eqnarray}\label{eq: I_1}
I_1&=&\int \grad u_{\epsilon}(t,x)\cdot a\grad\left(\frac{\phi_{R}(x)}{u_{\epsilon}(t,x)}\right)dx\nonumber\\
&=&\int \grad \log u_{\epsilon}(t,x)\cdot a\grad \phi_{R}(x)dx-\int \grad \log u_{\epsilon}(t,x)\cdot a\grad \log u_{\epsilon}(t,x) \phi_{R}(x)dx.
\end{eqnarray}
Using the ellipticity condition, we obtain the following
\begin{eqnarray*}
0&\leq&\int(\sqrt{\phi_{R}(x)}\grad \log u_{\epsilon}(t,x)-\frac{\grad \phi_{R}(x)}{\sqrt{\phi_{R}(x)}})\cdot a(\sqrt{\phi_{R}(x)}\grad \log u_{\epsilon}(t,x)-\frac{\grad \phi_{R}(x)}{\sqrt{\phi_{R}(x)}})\,dx\\
&=&\int \grad \phi_R(x)\cdot a\grad \phi_R(x) \phi_R(x)^{-1}dx-2\int \grad \log u_{\epsilon}(t,x)\cdot a \grad \phi_R(x) dx\\
&+&\int \phi_R(x)\grad \log u_{\epsilon}(t,x)\cdot a \grad \log u_{\epsilon}(t,x) dx.
\end{eqnarray*}
Rearranging the above and using the ellipticity condition again, we obtain 
\begin{equation*}
\int \grad \log u_{\epsilon}(t,x)\cdot a \grad \phi_R(x) dx \leq c_4|B_R(0)|R^{-2}+c_5\int|\grad \log u_{\epsilon}(t,x)|^2\phi_{R}(x)dx.
\end{equation*}
To obtain the above inequality, we have also used the following 
\begin{equation*}
\grad \phi_R(x)=2\frac{\sqrt{\phi_R(x)}}{R|x|}x,
\end{equation*}
with $|\phi_R(x)|\leq 1$. We now use the ellipticity condition once more and the above to bound $I_1$ as follows:
\begin{equation*}
I_1\leq c_4|B_R(0)|R^{-2}-c_6\int|\grad \log u_{\epsilon}(t,x)|^2\phi_{R}(x)dx.
\end{equation*}
See \cite{SS} where similar arguments are used.
Now using (\ref{poincare}) and the fact that $\mu(R)\asymp B_R(0)$, we obtain
\begin{equation*}
-R^2I_1\geq -c_7|B_R(0)|+c_8\int\phi_{R}(x)|\log u_{\epsilon}(t,x)-G_\epsilon(t)|^2dx.
\end{equation*}
Here $\mu(R)\asymp B_R(0)$ means that there is a constant $c>0$ such that $c^{-1}\mu(R)\leq B_R(0)\leq c\mu (R)$.  So combining the above, inequality (\ref{G1}) reduces to
\begin{equation}\label{Ma2}
G'_\epsilon(t)\geq -c_9+c_{10}\mu(R)^{-1}\int|\log u_{\epsilon}(t,x)-G_\epsilon(t)|^2\phi_{R}(x)dx.
\end{equation}
Let $D_t=\{x\in B_{R/2}(0): u_{\epsilon}(t,x)\geq e^{-K}\}$, where $K$ is a positive constant to be chosen later.  
By choosing $t_1$ small and using Proposition \ref{generaltightness}, we obtain
\begin{equation}\label{Ma1}
\bP^x(\sup_{s\leq t_1R^2}|X_s-X_0|\geq R/2)\leq \frac{1}{2}.
\end{equation}
 Using (\ref{Ma1}), we obtain
\begin{eqnarray*}
\int_{B(0,R/2)}p(tR^2,z,y)dy&\geq&1-\bP^z(\sup_{s\leq t_1R^2}|X_s-X_0|\geq R/2)\\
&\geq&1-1/2=1/2.\hskip10mm {\rm for}\hskip3mm t\leq t_1.
\end{eqnarray*}
So for $t\leq t_1$, we have 
\begin{eqnarray*}
\frac{|B_{R/2}(0)|}{2}&\leq& \int_{B_{R/2}(0)}u(t,x)dx\\
&=&\int_{D_t}u(t,x)dx+\int_{B_{R/2}(0)-D_t}u(t,x)dx\\
&\leq&\int_{D_t}u(t,x)dx+\int_{B_{R/2}(0)-D_t}u_{\epsilon}(t,x)dx.
\end{eqnarray*}
Note that on $B_{R/2}(0)-D_t$, we have $u_\epsilon(t,x)\leq e^{-K}$.  Choosing $K$ such that $e^{-K}=1/4$ and using $p(tR^2,x,y)\leq c_{11}t^{-d/2}R^{-d}$. This upper bound can be obtain by using an argument very similar to that of the proof of Proposition \ref{generalupper}  , we obtain
\begin{equation*}
\frac{|B_{R/2}(0)}{2}\leq c_{12}|D_t|t^{-d/2}+\frac{|B_{R/2}(0)|}{4}.
\end{equation*}
We thus obtain 
\begin{equation*}
|D_t|\geq \frac{t^{d/2}|B_{R/2}(0)|}{c_{13}}\quad{\rm for}\quad t\leq t_1.
\end{equation*}
Recall that 
\begin{equation*}
G_{\epsilon}(t)=\mu(R)^{-1}\int (\log u_\epsilon(t,x))\phi_R(x)\,dx.
\end{equation*}
Note that since $\epsilon$ is small, we can assume that it satisfies $\epsilon \leq c_{14}t^{-d/2}R^{-d}$. So for $t\leq t_1$, we can use the bound  $p(tR^2,x,y)\leq c_{15}t^{-d/2}R^{-d}$ and the above inequality to conclude that $G_\epsilon(t)$ is bounded above by a constant which we denote by $\bar{G}$. 

Since on $D_t$, $\log u_{\epsilon}(t,x)\geq -K$, we have only four possibilities:
\begin{itemize}
\item[(a)] If $\log u_{\epsilon}(t,x)>0$ and $G_\epsilon(t)\leq 0$, then $(\log u_{\epsilon}(t,x)-G_\epsilon(t))^2\geq G_\epsilon(t)^2$.\\
\item[(b)] If $\log u_{\epsilon}(t,x)>0$ and $0<G_\epsilon(t)\leq \bar{G}$, then $(\log u_{\epsilon}(t,x)-G_\epsilon(t))^2\geq 0\geq G_\epsilon(t)^2-\bar{G}^2$.\\
\item[(c)] If $-K\leq\log u_{\epsilon}(t,x)\leq 0$ and $|G_\epsilon(t)|<2K$, then $(\log u_{\epsilon}(t,x)-G_\epsilon(t))^2\geq 0\geq \frac{G_\epsilon(t)^2}{4}-K^2$.\\
\item[(d)] If $-K\leq\log u_{\epsilon}(t,x)\leq 0$ and $|G_\epsilon(t)|\geq2K$, then $(\log u_{\epsilon}(t,x)-G_\epsilon(t))^2\geq \frac{G_\epsilon(t)^2}{4}$.
\end{itemize}
We can therefore conclude that there exist positive constants $c_{16}$ and $c_{17}$ such that on $D_t$, 
\begin{equation*}
(\log u_\epsilon(t,x)-G_\epsilon(t))^2\geq -c_{16}+c_{17}G_\epsilon(t)^2.
\end{equation*}
Using the above and the fact that $\mu(R)\asymp B_{R/2}(0)$, inequality (\ref{Ma2}) then reduces to 
\begin{equation}\label{Ma4}
G'_\epsilon(t)\geq -c_{18}+c_{19}G_\epsilon(t)^2\,\,\,{\rm for\,\,\, all}\,\,\, t\in[t_1/2,t_1],
\end{equation}
where $c_{18}$ and $c_{19}$ are independent of $\epsilon$. Also note that $t_1$ is small and can be taken  to be less than one.  See the proof of Proposition 4.9 of \cite{BBCK} or the proof of Theorem 3.4 of \cite{CKK} for details.  Assume that $G_\epsilon(t_1)\leq -c_{18}-2(c_{18}/c_{19})^{1/2}$. We can now write
\begin{equation}\label{Ma5}
G'_\epsilon(t)\geq -c_{18}
\end{equation}
and use some calculus to show that $(\ref{Ma4})$ reduces to $G'_\epsilon(t)\geq (3/4)c_{19}G_\epsilon(t)^2$ and that $G_\epsilon({t_1}/2)<0$. This in turn implies that $G_\epsilon(t_1)\geq -8/(3c_{19})$.  We have thus obtain 
\begin{equation}
G_\epsilon(t_1)\geq -c_{20}\hskip10mm {\rm where}\,\,\,c_{20}=\sup\{c_{18}+2(c_{18}/c_{19})^{1/2},8/3c_{19}\}.
\end{equation}
Choose $y\in B_R(0)$.  By the semigroup property, we have 
\begin{eqnarray*}
|B_R(0)|p(2t_1R^2,z,y)&=&|B_R(0)|^{-1}\int|B_R(0)|p(t_1R^2,z,x)|B_R(0)|p(t_1R^2,x,y)dx\\
&\geq&|B_R(0)|^{-1}\int|B_R(0)|p(t_1R^2,z,x)|B_R(0)|p(t_1R^2,x,y)\phi_{R}(x)dx.
\end{eqnarray*}
Applying logarithm to the above and using Jensen's inequality we obtain
\begin{eqnarray*}
\log(|B_R(0)|p(2t_1R^2,z,y))&\geq&\log(\mu(R)|B_R(0)|^{-1})\\
&+&\mu(R)^{-1}\int\log[|B_R(0)|p(t_1R^2,z,x)]\phi_{ R}(x)dx\\
&+&\mu(R)^{-1}\int\log[|B_R(0)|p(t_1R^2,x,y)]\phi_{ R}(x)dx.
\end{eqnarray*}
Using the fact that $G_\epsilon(t_1)\geq -c_{20}$ after taking the limit $\epsilon \rightarrow 0$ as in the proof of Lemma 3.3.3 of \cite{Dav}, the above reduces to $\log[|B_R(0)|p(2t_1R^2,z,y)]\geq -c_{21}$. Set $t'=t_1R^2$. We have thus obtain $p(2t',z,y)\geq c_1{t'}^{-d/2}$ for $|z-y|^2\leq 2\theta t'$ which is the desired result with $\theta=2/{t_1}$.
\qed

\section{Some estimates}

The following estimates will be crucial for the proof of the regularity theorem and the Harnack inequality.

\begin{proposition}\label{Hit1}
Let $x_0\in \bR^d\backslash \sN$. Then the following holds
 \begin{enumerate}[(a)]
\item There exist constants $c_1$, $c_2$  and $r_0 $ such that $\bE^x\tau_{B(x_0,r)}\leq c_1r^2$, for $x\in B(x_0,r)$ and $r>0$ and $\bE^x\tau_{B(x_0,r)}\geq c_2r^2$ for $x\in B(x_0,r/2)$ and $r\in(0,r_0]$.
\item for any $A\subset B(x_0,3r/4)$, there exists some positive constant $c_3$ such that $\begin{displaystyle}\bP^x(T_A<\tau_{B(x_0,r)})\geq ~\frac{c_3|A|}{r^d}\end{displaystyle}$ for $x\in B(x_0,r/2)$ and $r\in(0,r_1]$ where $r_1$ is some positive constant.
\end{enumerate}
\end{proposition}

\proof{ Let $C\subset B(x_0,2r)\backslash B(x_0,r)$.  We can then write
\begin{eqnarray*}
\bP^x(\tau_{B(x_0,r)}<t)&\geq&\bP^x(\tau_{B(x_0,r)}<t, X_t\in C)\\
&\geq& \int_{C}p(t,x,y)dy\\
&\geq&c_4|C|t^{-d/2},
\end{eqnarray*}
where we used Theorem \ref{theo2} in the last inequality. Taking $|C|=c_5r^d$  and  $t=c_6r^2$ we obtain upon choosing $c_6=(2c_4c_5)^{2/d}$,
\[\bP^x(\tau_{B(x_0,r)}<c_6r^2)\geq \frac{1}{2}.\]
Let $m$ be a positive integer. By the Markov property and using induction, we obtain
\[\bP^x(\tau_{B(x_0,r)}\geq mc_6r^2)\leq 2^{-m}.\]
We can now obtain $\bE^x\tau_{B(x_0,r)}\leq c_2r^2$ from the above. Let $t=c_7r^2$, then by Proposition $\ref{generaltightness}$, we have $\bP^x(\tau_{B(x_0,r)}\leq t)\leq \bP^x(\tau_{B(x,r/2)}\leq t)\leq 1/2$ for $c_7$ small enough.  We thus have 
\begin{eqnarray*}
\bE^x\tau_{B(x_0,r)}&\geq&t\bP^x(\tau_{B(x_0,r)}\geq t)\\
&\geq&c_8r^2.
\end{eqnarray*}
For part(b), since we need to prove a lower bound, it suffices to obtain the result for small jumps(less than $\lambda$) only.  The more general result follows from the following fact:
\begin{eqnarray*}
\bP^x(T_A<\tau_{B(x_0,r)})&=&\bP^x(T_A<\tau_{B(x_0,r)}, U_1\leq t)+\bP^x(T_A<\tau_{B(x_0,r)}, U_1>t)\\
&\geq&\bP^x(T_A<\tau_{B(x_0,r)}\leq t, U_1> t)\\
&\geq&e^{-(\sup N) t}\bP^x(T^\lambda_A<\tau^\lambda_{B(x_0,r)}),
\end{eqnarray*}
where $U_1$  and $N(x)$ are defined in Remark \ref{remark2}. The stopping times $T^\lambda_A$ and $\tau^\lambda_{B(x_0,r)}$ are defined in a similar way to $T_A$ and $\tau_{B(x_0,r)}$ respectively but for processes with jumps less than $\lambda$.  So from now on, we assume that $X$ has jumps less than $\lambda$. For $t$ fixed, we can now write
\begin{eqnarray}\label{Heq1}
\bP^x(T_A<\tau_{B(x_0,r)})&\geq&\int_Ap(t,x,y)dy\nonumber\\
&-&\int_A\bE^x[\bE^{X_{\tau_{B_r(x_0)}}}[p(t-\tau_{B_r(x_0)}, X_{\tau_{B_r(x_0)}},y)]1_{(\tau_{B_r(x_0)}<t)}]dy.
\end{eqnarray}
Since our process is assumed to have small jumps only, we can use (\ref{GT1}) to obtain, for $t$ sufficiently small,
\begin{equation}\label{Heq2}
p(t-\tau_{B_r(x_0)}, X_{\tau_{B_r(x_0)}},y)\leq c_9t^\gamma e^{-c_{10}|X_{\tau_{B_r(x_0)}}-y|},
\end{equation}
where we have taken $\lambda$ small enough so that $\gamma=|X_{\tau_{B_r(x_0)}}-y|/12\lambda-d/2>0$ whenever $|X_{\tau_{B_r(x_0)}}-y|>r/4$.  We now use the lower bound given by Theorem \ref{theo2} to reduce $(\ref{Heq1})$ to  
\begin{equation}
\bP^x(T_A<\tau_{B(x_0,r)})\geq \frac{c_{11}|A|}{r^{d}}.
\end{equation}
We have taken $t=c_{12}r^2$, where $r\in(0,r_1]$ and $r_1$ is a small constant so that the right hand side of (\ref{Heq2}) is less than a positive fraction of the lower bound on the heat kernel $p(t,x,y)$.

\qed}

\section{The Regularity Theorem}

We will need the following L\'evy system formula for our process $X$.  The proof is the same as that of Lemma~4.7 in \cite{CK}. So we omit it here.
\begin{proposition}\label{levysystem}
If $A$ and $B$ are disjoint Borel sets, then for each $x\in \bR^d\backslash \sN,$
\begin{equation}\label{levy1}
\bE^x\sum_{s\leq t}1_{(X_{s-}\in A,\,X_s\in B)}=\bE^x\int_0^t\int_B1_A(X_s)J(X_s,y)dyds.
\end{equation}
\end{proposition}

The proof of the following is based on the proof of the regularity theorem in \cite{BL1}.  For the sake of completeness, we give a proof here.

\longproof{of Theorem \ref{theo3}}: 
Let us suppose $u$ is bounded by $M$ in $\bR^d$ and $z_1\in B(z_0,R/2)\backslash \sN$. Set
\[ s_n=\theta_1a^n, \hskip 10mm  r_n= \theta_2 \rho^{n}, \hskip10mm {\rm for\,\,\,\,}n\in \bN,\]
where $a < 1$, $\rho<\frac{1}{2}$, and $\theta_1 \ge 2M$ are constants to be chosen later. We choose $\theta_2$ small enough that $B(z_1, 2r_1)\subset B(z_0,R/2)$. Write $B_n=B(z_1,r_n)$ and $\tau_n=\tau_{B_n}$. Set
\[M_n=\sup_{x\in B_n} u(x), \hskip 10mm m_n=\inf_{x\in B_n} u(x).\]
H\"older continuity will follow from the fact that $M_n-m_n \leq s_n$ for all $n$ which will be proved by induction.  Let $n_0$ be a positive number to chosen later and suppose $M_i-m_i\leq s_i$ for all $i=1,2,...,n$, where $n\geq n_0$; we want to show 
\[M_{n+1}-m_{n+1} \leq s_{n+1}.\]
\
Let
\[ A_n=\{ z \in B_{n}: u(z)\leq (M_n+m_n)/2\}.\]
We may suppose that $|A_n|/|B_{n}|\geq \frac{1}{2}$, for if not, we can look at $M-u$ instead. Let $A$ be compact subset of $A_n$ such that $|A|/|B_n|\geq 1/3$. By  Proposition \ref{Hit1}, there exists $c_1$ such that
\begin{equation}\label{R0}
\bP^x(T_A\leq \tau_{n})\geq c_1
\end{equation}
 for all $x\in B_{n+1}$.  Let $\epsilon>0$ and pick $y\,,z \in B_{n+1}$ such that $u(y)\leq m_{n+1}+\epsilon$ and $u(z)\geq M_{n+1}-\epsilon$.  Since $\epsilon>0$ is arbitrary, showing that $u(z)-u(y)\leq s_{n+1}$ will imply $M_{n+1}-m_{n+1}\leq s_{n+1}$ as desired. 

By optional stopping,
\begin{eqnarray}\label{R1}
u(z)-u(y)&=&\bE^z[u(X_{T_A})-u(y); T_A \le \tau_n] \nonumber\\
&+& \bE^z[u(X_{\tau_n})-u(y); \tau_n \le T_A, X_{\tau_n}\in B_{n-1}] \nonumber\\
&+& \sum_{i=1}^{n-2} \bE^z[u(X_{\tau_n})-u(y); \tau_n \le T_A, X_{\tau_n}\in B_{n-i-1}-B_{n-i}] \nonumber\\
&+& \bE^z[u(X_{\tau_n})-u(y);  \tau_n \le T_A, X_{\tau_n}\notin B_1]\nonumber\\
&=& I_1+I_2+I_3+I_4.
\end{eqnarray}
By optional stopping and the L\'evy system formula (\ref{levy1}),
\begin{eqnarray}\label{R2}
\sup_{y\in B_{n+1}}\bP^y(X_{\tau_n}\notin B_{n-i})&\leq& \sup_{y\in B_{n+1}}\bE^y\tau_n\int_{|x-y|>r_{n-i}-r_n}J(x,y)dx \nonumber\\
&=& \sup_{y\in B_{n+1}}\bE^y\tau_n[\int_{|x-y|>1}J(x,y)dx+\int_{r_{n-i}-r_n<|x-y|<1}J(x,y)dx]\nonumber\\
&\leq&c_2r_n^2+c_3\left(\frac{\rho^i}{1-\rho^i}\right)^2.
\end{eqnarray}
See the proof of Proposition 3.5 of \cite{BL1} where a similar argument is used.
We have also used Proposition \ref{Hit1}(a), Assumption \ref{ass2}(a) and the fact that $1<\frac{|x-y|}{r_{n-i}-r_n}$ in the above computations.
The first term on the right of (\ref{R1}) is bounded by 
\begin{equation}\label{R3}
\left(\frac{M_n+m_n}{2}-m_n\right)\bP^y(T_A \le \tau_n) \leq \frac{1}{2}s_n\bP^y(T_A\le \tau_n).\\
\end{equation}
The second term is bounded by
\begin{equation}\label{R4}
(M_{n-1}-m_{n-1})\bP^y(\tau_n \le T_A)=(M_{n-1}-m_{n-1})(1-\bP^y(T_A\le \tau_n)) \leq s_{n-1}(1-\bP^y(T_A \leq \tau_n)).\\
\end{equation}
Let $\rho^2=\frac{\sqrt{a}}{2}\wedge \sqrt{\frac{3c_1a}{128c_3}}$.  Using (\ref{R2}), the third term is bounded by
\begin{eqnarray*}
 \sum_{i=1}^{n-2}(M_{n-i-1}&-&m_{n-i-1})\bP^y(X_{\tau_n}\notin B_{n-i}) \nonumber\\
 &\leq& c_2\sum_{i=1}^{n-2}s_{n-i-1}r_n^2+c_3\sum_{i=1}^{n-2}s_{n-i-1}\rho^{2i}\nonumber\\
&\leq& s_{n-1}[\frac{c_2a^2\theta_2^2\rho^{2n}}{a^n(1-a)}+\frac{c_3\rho^2/a}{1-\rho^2/a}]\nonumber\\
&\leq&s_{n-1}[\frac{a^2\theta_2^2c_4}{1-a}+\frac{c_1}{32}].
 \end{eqnarray*}
where we can take $n_0$ bigger if necessary so that the last inequality holds. We also choose $\theta_2$ smaller if necessary so that $\theta_2\leq \frac{1}{4}\sqrt{\frac{c_1(1-a)}{2a^2c_4}}$ and obtain
\begin{equation}\label{R5}
I_3 \leq \frac{s_{n-1}c_1}{16}
\end{equation}
The fourth term can be bounded similarly 
\begin{eqnarray*}
 I_4 \leq 2M \bP^y(X_{\tau_n} \notin B_1)&\leq&2M [c_2r_n^2+c_3\rho^{2(n-1)}]\nonumber\\
 &\leq&\theta_1[c_2a^{4n}\theta_2^2+c_3a^{4n-4}].
 \end{eqnarray*}
 By choosing $n_0$ bigger if necessary, we have, for $n\geq n_0$, the above yields
 \begin{equation}\label{R6}
 I_4\leq \frac{s_{n-1}c_1}{8}.
 \end{equation}
Inequalities (\ref{R3})-(\ref{R6}) give the following:
 \begin{eqnarray*}
 u(y)-u(z)&\leq& \frac{1}{2}as_{n-1}\bP^y(T_A\le \tau_n)+ s_{n-1}(1-\bP^y(T_A \leq \tau_n)) +s_{n-1}(\frac{c_1}{16}+\frac{c_1}{8}) .\nonumber
 \end{eqnarray*}
 Using the fact that $a$ is less than one, we obtain
 \begin{eqnarray*}
 u(z)-u(y)&\leq& \frac{s_n}{a}\Big[1-\frac{\bP^y(T_A<\tau_n)}{2}+\frac{c_1}{16}+\frac{c_1}{8}\Big]\\
 &\leq&\frac{s_n}{a}[1-\frac{5c_1}{16}].\\
 \end{eqnarray*} 
Now let us pick $a$ as follows:
 \[  a=\sqrt{1-\frac{5c_1}{16}}.\]
This yields
\begin{equation*}
u(z)-u(y)\leq s_na=s_{n+1}.
\end{equation*}
\qed


\section{The Harnack Inequality}

We start this section with the following proposition which will be used in the proof of the Harnack inequality.
\begin{proposition}\label{Har1}
Let $x_0\in \bR^d$ and $r\leq r_0$, where $r_0$ is a positive constant. Then there exists $c_1$ depending on $\kappa,\, K_is$ and $\Lambda$ such that if $z\in B(x_0,r/4)$ and $H$ is a bounded non-negative function supported in $B(x_0,r)^c$, then
\[\bE^{x_0}H(X_{\tau_{B(x_0,r/2)}})\leq c_1k_r\bE^zH(X_{\tau_{B(x_0,r/2)}}).\]
\end{proposition}
\proof{ By linearity and a limit argument, it suffices to show to consider only $H(x)=1_C(x)$ for a set $C$ contained in $B(x_0,r)^c$. From Assumption \ref{ass2}(c), we have $J(w,v)\leq k_rJ(y,w)$ for all $w,\,y\in B(x_0, r/2)$ and $v\in B(x_0,r)^c$.  Hence, we have
\begin{equation}\label{a1}
\sup_{y\in B(x_0,r/2)}J(y,v)\leq k_r \inf_{y\in B(x_0,r/2)}J(y,v).
\end{equation}
By optional stopping and the L\'evy system formula, we have 
\begin{eqnarray*}
\bE^z1_{(X_{t\wedge \tau_{B(x_0,r/2)}}\in C)}&=&\bE^z\sum_{s\leq{t\wedge \tau_{B(x_0,\frac{r}{2})}}}1_{(|X_s-X_{s-}|\geq \frac{r}{2}, X_s\in C)}.\\
&=&\bE^z\int_0^{t\wedge \tau_{B(x_0,\frac{r}{2})}}\int_C J(X_s,v)dvds.\\
&\geq&\bE^z(t\wedge \tau_{B(x_0,\frac{r}{2})})\int_C \inf_{y\in B(x_0,\frac{r}{2})}J(y,v)dv.
\end{eqnarray*}
Letting $t\to\infty$ and using the dominated convergence theorem on the left and monotone convergence on the right, we obtain
\[\bP^z(X_{\tau_{B(x_0,\frac{r}{2})}}\in C)\geq \bE^z\tau_{B(x_0,\frac{r}{2})}\int_C\inf_{y\in{B(x_0,\frac{r}{2})}}J(y,v)dv.\]
Since $ \bE^z\tau_{B(x_0,\frac{r}{2})}\geq \bE^z\tau_{B(z,\frac{r}{4})}$, we have 
\begin{equation}\label{a2}
\bP^z(X_{\tau_{B(x_0,\frac{r}{2})}}\in C)\geq \bE^z\tau_{B(z,\frac{r}{4})}\int_C\inf_{y\in{B(x_0,\frac{r}{2})}}J(y,v)dv. 
\end{equation}
Similarly we have 
\begin{equation}\label{a3}
\bP^{x_0}(X_{\tau_{B(x_0,\frac{r}{2})}}\in C)\leq \bE^{x_0}\tau_{B(x_0,\frac{r}{2})}\int_C\sup_{y\in{B(x_0,\frac{r}{2})}}J(y,v)dv.
\end{equation}
Combining inequalities (\ref{a1}), (\ref{a2}) and (\ref{a3}) and using Proposition \ref{Hit1}(a), we get our result.
\qed}

\longproof{of theorem \ref{theo4}}:
By looking at $u+\epsilon$ and letting $\epsilon \downarrow 0$, we may suppose that $u$ is bounded below by a positive constant. Also, by looking at $au$, for a suitable $a$, we may suppose that $\inf_{B(z_0,R/2)}u \in [1/4,1]$.
 \
 We want to bound $u$ above in $B(z_0,R/2)$ by a constant not depending on $u$.  Our proof is by contradiction.
 
 Since $u$ is continuous, we can choose $z_1\in B(z_0,R/2)$ such that $u(z_1)=\frac{1}{3}$.  Let $r_i=r_1i^{-2}$ where $r_1<r_0$ is a chosen constant so that $\sum_{i=1}r_i< R/8.$ Recall that from  Proposition \ref{Har1}, there exists $c_1$ such that if $r<r_0$, $z\in B(x,r/4)$ and $H$ is a bounded non-negative function supported in $B(x,r)^c$, then
\begin{equation}\label{H1}
\bE^xH(X_{\tau_{B(x,r/2)}}) \leq c_1k_r\bE^zH(X_{\tau_{B(x,r/2)}}).
\end{equation}
We will also use Proposition \ref{Hit1}(b) which says that if $A\subset B(x,3r/4)$, then there exists a constant $c_2$ such that 
\begin{equation}\label{H2}
\bP^x(T_A<\tau_{B(x,r)})\geq \frac{c_2|A|}{r^d}.
\end{equation}
Let  $\eta$ be a constant to be chosen later and let $\xi$ be defined as follows
\[\xi= \frac{1}{2}\wedge c_1^{-1}\eta.\]
Let $c_3$ be a positive constant to be chosen later. Once this constant has been chosen, we suppose that there exists $x_1\in B(z_0, R/2)$ with $u(x_1)=L_1$ for some $L_1$ large enough so that we have
\begin{equation}\label{H5}
\frac{c_2\xi L_1e^{c_3j}r_j^{d+\beta}}{\kappa2^{2d+1}}>2,\quad{\rm for\,\,all}\quad j.
\end{equation}
The constants $\beta$ and $\kappa$ are from Assumption \ref{ass2}(c).

We will show that there exists a sequence $\{(x_j,L_j)\}$ with $x_{j+1}\in \overline{B(x_j, r_j)}\subset B(x_j,2r_j)\subset B(z_0,3R/4)$ with
\begin{equation*}
L_j=u(x_j)\quad{\rm and}\quad L_j\geq L_1e^{c_3j}.
\end{equation*}
This would imply that $L_j\rightarrow \infty$ as $j\rightarrow \infty$ contradicting the fact that $u$ is bounded.  Suppose that we already have $x_1,x_2,...,x_i$ such that the above condition is satisfied.  We will show that there exists $x_{i+1}\in \overline{B(x_i,r_i)}\in B(x_i,2r_i)$ such that $L_{i+1}=u(x_{i+1})$ and $L_{i+1}\geq L_1e^{c_3(i+1)}$. Define 
\begin{equation*}
A=\{ y\in B(x_i,r_i/4);  u(y)\geq \frac{\xi L_ir_i^\beta}{\kappa} \}. 
\end{equation*}

 We are going to show that $\begin{displaystyle} |A|\leq \frac{1}{2}|B(x_i, r_i/4)| \end{displaystyle}
$.  To prove this fact, we suppose the contrary. Choose a compact set $A' \subset A$ with $\begin{displaystyle} |A'| > \frac{1}{2}|B(x_i,r_i/4)| \end{displaystyle}$.

By optional stopping, (\ref{H2}), the induction hypotheses and the fact that $R<1$,
\begin{eqnarray*}
\frac{1}{3}&=& u(z_1)\geq \bE^{z_1}[u(X_{T_{A'}\wedge \tau_{B(z_0,R)}}); T_{A'} < \tau_{B(z_0,R)}]\\
&\geq& \frac{\xi L_ir_i^\beta}{\kappa}\bP^{z_1}(T_{A'}<\tau_{B(z_0,R)})\\
&\geq& \frac{\xi L_ir_i^\beta}{\kappa} \frac{c_2|A'|}{R^d}\\
&=& \frac{c_2\xi L_ir_i^\beta}{\kappa}\frac{|A'|}{(r_i/4)^d}\frac{(r_i/4)^d}{R^d}\\
&\geq& 2.
\end{eqnarray*}
This is a contradiction.  Therefore $\begin{displaystyle}|A|\leq \frac{1}{2} |B(x_i, r_i/4)|\end{displaystyle}$.
So we can find a compact set $E$ such that $E \subset B(x_i, r_i/4)-A$ and $|E|\geq \frac{1}{3}|B(x_i,r_i/4)|$. Let us write $\tau_{r_i}$ for $\tau_{B(x_i,r_i/2)}$.  From (\ref{H2}) we have $\bP^{x_i}(T_E<\tau_{r_i})\geq c_4$ where $c_4$ is some positive constant.
Let $M=\sup_{x\in B(x_i,r_i)}u(x)$.  We then have
\begin{eqnarray}\label{H3}
L_i=u(x_i)&=&\bE^{x_i}[u(X_{T_E\wedge \tau_{r_i}}); T_E < \tau_{r_i}]\nonumber \\
&+& \bE^{x_i}[u(X_{T_E\wedge \tau_{r_i}}); T_E > \tau_{r_i}, X_{\tau_{r_i}}\in B(x_i,r_i)]\nonumber \\
&+&\bE^{x_i}[u(X_{T_E\wedge \tau_{r_i}}); T_E > \tau_{r_i}, X_{\tau_{r_i}}\notin B(x_i,r_i)]\nonumber \\
&=& I_1+I_2+I_3.
\end{eqnarray}
Writing $p_i=\bP^{x_i}(T_E<\tau_{r_i})$, we see that the first two  terms are easily bounded as follows:
\[ I_1 \leq \frac{\xi L_ip_ir_i^\beta}{\kappa}, \hskip 15mm {\rm and}\hskip 15mmI_2\leq M(1-p_i).\]
To bound the third term, we prove $\bE^{x_i}[u(X_{\tau_{r_i}}); X_{\tau_{r_i}} \notin B(x_i,r_i)] \leq \eta L_i$.  If not, then by using (\ref{H1}), we will have, for all $y\in B(x_i,r_i/4)$, 
\begin{eqnarray*}
u(y)&\geq& \bE^yu(X_{\tau_{r_i}})\geq \bE^y[u(X_{\tau_{r_i}}); X_{\tau_{r_i}}\notin B(x_i,r_i)]\\
&\geq& \frac{1}{k_{r_i}c_1}\bE^{x_i}[u(X_{\tau_{r_i}});X_{\tau_{r_i}}\notin B(x_i,r_i)]>\frac{\eta L_i}{k_{r_i}c_1}>\frac{\xi L_ir_i^\beta}{\kappa}.
\end{eqnarray*}
contradicting the fact that $A\leq \frac{1}{2}|B(x_i,r_i/4)|$. Hence
\[ I_3 \leq \eta L_i.\]
So (\ref{H3}) becomes
\[L_i\leq \frac{\xi L_ip_ir_i^\beta}{\kappa}+M(1-p_i)+\eta L_i,\]
or 
\begin{equation}\label{H4}
\frac{M}{L_i}\geq \frac{1-\eta-\xi p_ir_i^\beta/\kappa}{1-p_i}.
\end{equation}
Choosing $\eta=\frac{c_4}{4}$ and using the definition of $\xi$ together with the fact that $p_i\geq c_5$ and $r_i^\beta/\kappa<1$, we see that there exists a constant $\gamma$ bounded below by a positive constant, such that the inequality ($\ref{H4}$) reduces to $M\geq L_i(1+\gamma)$. Therefore, there exists $x_{i+1}\in \overline{B(x_i,r_i)}$ with $u(x_{i+1})\geq L_i(1+\gamma)$. Setting $L_{i+1}=u(x_{i+1})$, we see that 
\begin{eqnarray*}
L_{i+1}&\geq& L_i(1+\gamma)\\
&=&L_ie^{\log(1+\gamma)}.
\end{eqnarray*}
The induction hypotheses is thus  satisfied by taking $c_3=\log(1+\gamma)$.
\qed

\section*{Acknowledgments}
The author wishes to thank Prof. Richard F. Bass for suggesting this problem and all his encouragement during the preparation of this paper. The author also thanks Prof. Toshihiro Uemura and an anonymous referee for their helpful comments.
\begin{small}
\end{small}


\end{document}